\documentclass[11pt]{amsart}
\usepackage{latexsym}
\usepackage{hyperref}
                \usepackage{amscd}
                \usepackage{amssymb}
                \usepackage{amsmath}
                \usepackage{amsthm}
                \usepackage{enumerate}
                \usepackage{verbatim}

\newtheorem{thm}{Theorem}[subsection]

\newtheorem{Example}[thm]{Example}

\newtheorem{Counterexample}[thm]{Counterexample}

\newtheorem{remark}[thm]{Remark}

\newtheorem{Fact}[thm]{Fact}

\newtheorem{Nothing}[thm]{$\!\!\!$}

\begin{document}
\abovedisplayskip=6pt plus3pt minus3pt \belowdisplayskip=6pt
plus3pt minus3pt
\title[A note on nonnegative curvature and (rational) homotopy type]
{A note on nonnegative curvature and (rational) homotopy type}

\author{Anand Dessai and Wilderich Tuschmann}
\date{23 March 2011}

\address{Anand Dessai\\
Department of Mathematics\\
University of Fribourg\\
Chemin du Mus\'ee 23\\
CH-1700 Fribourg, Switzerland}
\email{anand.dessai@unifr.ch}

\address{Wilderich Tuschmann\\
Department of Mathematics\\
Karlsruhe Institute of Technology\\
Kaiserstrasse 89-93\\
D-76133 Karlsruhe, Germany}
\email{tuschmann@kit.edu}

\maketitle

\section{Introduction}

\noindent The study of Riemannian manifolds with given curvature
bounds and their topology is a central theme of global differential
geometry, and much attention has here been paid to the interplay
between curvature and homotopy type. In this regard there is also the following
question, which has been raised by F\'elix, Oprea, and Tanr\'e (see
Question~6.17 in their book \cite{FOT08}):

\

\noindent {\bf Question 1 ([FOT08], p.~249)} {\sl If two closed
manifolds have the same (rational) homotopy type and one manifold
has a metric of non-negative sectional curvature, does the other?}

\

In turns out that the answer to the above question is negative for
many classes of manifolds, even if one restricts to positive
sectional curvature and the simply connected case. The main purpose
of the present note - further details will appear elsewhere in a broader context - consists in providing some of these examples and
in proposing a substitute for Question~1.

\

Concerning the first point, examples of closed simply connected
manifolds $M$ with positive sectional curvature that have homotopy
equivalent cousins $N$ which do not admit a metric of nonnegative
scalar curvature (and hence nonnegative sectional curvature) can be
obtained in the following way:

\

\noindent 1. Let $M=S^n$ be the round sphere in dimension $n=8k+1$
or $8k+2$ and let $N$ be a homotopy $n$-sphere with non-vanishing
$\alpha $-invariant. Then $N$ does not admit a metric of positive
scalar curvature by the work of Hitchin \cite{Hi74}, and the
well-known deformation properties of scalar curvature imply that any
metric of nonnegative sectional curvature on $N$ must be flat.
However, since homotopy spheres do not carry flat metrics, $N$ does
not admit a metric of nonnegative sectional curvature.

Actually, the homotopy sphere $N$ does not even admit a metric of
nonnegative scalar curvature: As before, any such metric $g$ must be
scalar-flat. The non-vanishing of the $\alpha$-invariant now implies that $N$ carries a non-trivial parallel spinor and
 that $(N,g)$ must have special holonomy (compare \cite{Hi74}).
However, since homotopy spheres do always have generic holonomy, $N$
doesn't carry a metric of nonnegative scalar curvature.

Essentially the same reasoning can be applied to any spin manifold
$M$ of positive sectional curvature in dimension $8k+1$ or $8k+2$ to
construct a manifold which is homotopy equivalent to $M$ but which
does not admit a metric of nonnegative sectional curvature.

\

\noindent 2. Let $M=\mathbb{H} P^m$ be the quaternionic projective
space of dimension $n=4m\geq 8$ equipped with its natural metric of
positive sectional curvature. Since the universal expressions for
the $\hat A$-genus and the signature in terms of Pontrjagin numbers
are linearly independent, we may, using surgery theory, choose a
manifold $N$ with non-vanishing $\hat A$-genus in the homotopy type
of $M$. Thus $M$ has positive sectional curvature whereas $N$ does
not admit any metric of positive scalar curvature. Moreover,
holonomy considerations as above show that $N$ does not even admit a
metric of nonnegative scalar curvature.

\

In search of suitable substitutes for Question~1, given the above
examples, one might first restrict Question~1 to manifolds which
already admit metrics of nonnegative or even positive scalar
curvature, e.g.: {\sl If two closed manifolds of positive scalar
curvature have the same (rational) homotopy type and one manifold
has a metric of non-negative sectional curvature, does the other?}

However, there is strong evidence that the answer to this question
might be negative as well. Indeed, if, for example, Stolz'
conjecture about the vanishing of the Witten genus of string
manifolds with positive Ricci curvature is true, one can construct
counterexamples.

On the other hand, one may also restrict Question~1 to manifolds
which a priori admit metrics with positive Ricci curvature:

\

\noindent {\bf Question 1$^\prime$} {\sl If two closed manifolds of
positive {\it Ricci} curvature have the same (rational) homotopy
type and one manifold has a metric of non-negative sectional
curvature, does the other?}

\

This question is at present completely open and, moreover, any kind of
answer to it would be interesting. In particular, a negative answer
will show that there exist manifolds of positive Ricci curvature
which do not admit metrics with nonnegative sectional curvature, but
whose homology satisfies the bounds given by Gromov's Betti number
theorem for nonnegatively curved manifolds.

\small
\bibliographystyle{alpha}

\end{document}